\documentclass[11pt]{amsart}

\usepackage{graphicx,amssymb,epstopdf,subcaption,mathtools,tikz-cd,pb-diagram,thmtools, thm-restate,pinlabel,pifont,enumitem,stmaryrd}

\newcommand{\nc}{\newcommand}
\nc{\dmo}{\DeclareMathOperator}
\nc{\nt}{\newtheorem}

\nt{theorem}{Theorem}
\nt{corollary}[theorem]{Corollary}
\nt{lemma}[theorem]{Lemma}
\nt{criterion}[theorem]{Criterion}
\nt*{theorem*}{Theorem}
\nt{question}[theorem]{Question}
\nt{problem}[question]{Problem}
\nt{conjecture}[question]{Conjecture}
\nt{proposition}[theorem]{Proposition}

\nc{\M}{\mathcal{M}}
\nc{\C}{\mathcal{C}}

\nc{\cut}{\!\ssearrow\!}

\dmo{\Diff}{Diff}
\dmo{\Mod}{Mod}
\dmo{\SMod}{SMod}
\dmo{\I}{\mathcal{I}}
\dmo{\SO}{SO}
\dmo{\Orth}{O}
\dmo{\Sp}{Sp}
\dmo{\SL}{SL}
\dmo{\GL}{GL}
\dmo{\PSp}{PSp}
\dmo{\PSL}{PSL}
\dmo{\Homeo}{Homeo}
\dmo{\Aut}{Aut}
\dmo{\Fix}{Fix}
\dmo{\SAut}{SAut}

\nc{\Z}{\mathbb Z}
\nc{\N}{\mathcal N}
\nc{\R}{\mathbb R}
\nc{\F}{\mathcal F}
\nc{\Sph}{\mathbb S}

\nc{\ga}{\gamma}
\nc{\de}{\delta}
\nc{\ep}{\epsilon}

\nc{\flm}{\lambda_{2}}

\nc{\normalclosure}[1]{\ensuremath{\left \langle \left \langle #1 \right \rangle \right \rangle}}

\nc{\margin}[1]{\marginpar{\scriptsize #1}}

\nc{\p}[1]{\bigskip\noindent\textbf{#1.}}
\nc{\lei}[1]{{\color{red} \sf  L: [#1]}}
\nc{\justin}[1]{{\color{green} \sf  J: [#1]}}

\title[Constraining mapping class group homomorphisms]{Constraining mapping class group homomorphisms using finite subgroups}
\author{Lei Chen}
\author{Justin Lanier}

\address{Lei Chen \\ Department of Mathematics\\ University of Maryland \\ 4176 Campus Drive\\ College Park, MD 20742 \\  chenlei@umd.edu}

\address{Justin Lanier \\ Department of Mathematics\\ University of Chicago \\ 5734 S. University Avenue \\ Chicago, IL 60637 \\  jlanier8@uchicago.edu}

\begin{document}

\maketitle

\vspace*{-4ex}

\begin{abstract}
We classify homomorphisms from mapping class groups by using finite subgroups. First, we give a new proof of a result of Aramayona--Souto that homomorphisms between mapping class groups of closed surfaces are trivial for a range of genera. Second, we show that only finitely many mapping class groups of closed surfaces have non-trivial homomorphisms into $\text{Homeo}(\mathbb{S}^n)$  for any $n$. We also prove that every homomorphism from $\text{Mod}(S_g)$ to $\text{Homeo}(\mathbb{S}^2)$ or $\text{Homeo}(\mathbb{S}^3)$ is trivial if $g\ge 3$, extending a result of Franks--Handel.
\end{abstract}


\vspace*{0in}

\section{Introduction}

Let $S_g$ be the connected, closed, orientable surface of genus $g$ and let $\Sph^n$ be the \mbox{$n$-sphere}. The groups $\Homeo^+(X)$ and $\Diff^+(X)$ are, respectively, the groups of orientation-preserving homeomorphisms and diffeomorphisms of a given orientable manifold $X$. The mapping class group $\Mod(X)$ is $\pi_0(\Homeo^+(X))$, the group of isotopy classes of orientation-preserving homeomorphisms of $X$.

In this paper, we analyze finite subgroups of $\Mod(S_g)$ to prove theorems about homomorphisms from $\Mod(S_g)$ to $\Mod(S_h)$ as well as to $\Homeo^+(\Sph^n)$ and $\Diff^+(\Sph^n)$. Much has been learned about homomorphisms of mapping class groups by analyzing torsion. Several examples appear later in the introduction, and other examples include work by Markovic, Mann--Wolff, and the first author \cite{Markovic,mann,ChenBraid}. For the problems we consider, it is not sufficient to analyze individual periodic elements in isolation. As a consequence, we proceed instead by analyzing non-cyclic finite subgroups of $\Mod(S_g)$. Our theorems are as follows.

\begin{theorem}
\label{thm:newAS}
For $g \geq 3$ and $0 \leq h< 2g-1$ with $h \neq g$, every homomorphism 
\[\phi: \Mod(S_g) \rightarrow \Mod(S_h)\] is trivial. When $g \geq 3$ is odd, the same conclusion holds for the bounds $0~\leq~h<~2g+1$ with $h \neq g$.
\end{theorem}
For $\Homeo^+(\Sph^n)$ as the target group, we prove the following theorem, suggested to us by Mattia Mecchia.
\begin{theorem}
\label{thm:finite}
For any $n$, there are only finitely many $g$ such that there exists a nontrivial homomorphism $\phi: \Mod(S_g)\to \Homeo^+(\mathbb{S}^n)$. 
\end{theorem}

We also effectivize Theorem~\ref{thm:finite} for some small values of $n$. Arguments treating some small values of $g$ were suggested to us by J\"{u}rgen M\"{u}ller.
\begin{theorem}
\label{thm:homeo}
For $g \geq 3$, every homomorphism from $\Mod(S_g)$ to $\Homeo^+(\Sph^2)$, $\Homeo^+(\R^3)$ or $\Homeo^+(\Sph^3)$ is trivial. For $g \geq 5$, every homomorphism from $\Mod(S_g)$ to $\Diff^+(\Sph^4)$ is trivial.
\end{theorem}

Note that since $\Mod(S_g)$ is perfect for $g \geq 3$, every homomorphism to a homeomorphism group has image in the subgroup of orientation-preserving homeomorphisms. Thus each of the above results also holds when the target group is the full homeomorphism group.

Since the action of $\Mod(S_g)$ on $H_1(S_g;\mathbb{Z})$ gives a surjective homomorphism $\Mod(S_g)\to \Sp(2g,\mathbb{Z})$, Theorem \ref{thm:finite} implies the following corollary.
\begin{corollary}
For any $n$, there are only finitely many $g$ such that there exists a nontrivial homomorphism $\phi: \Sp(2g,\mathbb{Z})\to \Homeo^+(\mathbb{S}^n)$.
\end{corollary} 
\noindent This result was originally proved by Zimmermann \cite{ZimmermannSp}, who additionally established the optimal bound on $g$.

Our proofs of Theorems \ref{thm:newAS}, \ref{thm:finite}, and \ref{thm:homeo} use the following result of the second author and Margalit. A group element is said to normally generate if its normal closure (the smallest normal subgroup containing this element) equals the group.

\begin{theorem}[Theorem~1.1, \cite{LM}]
\label{thm:LM}
For $g \geq 3$, every nontrivial periodic mapping class that is not a hyperelliptic involution normally generates $\Mod(S_g)$.
\end{theorem}

It follows immediately from this theorem that any homomorphism from $\Mod(S_g)$ that has a nontrivial nonhyperelliptic periodic element in its kernel is trivial. The strong constraints on homomorphisms provided by Theorem~\ref{thm:LM} are not, however, immediately applicable to the settings of Theorems~\ref{thm:newAS}, \ref{thm:finite}, and \ref{thm:homeo}. Under the hypotheses of each of these theorems, there are cases where, a priori, there could be a homomorphism to the target group where no nontrivial periodic mapping class of $\Mod(S_g)$ would lie in the kernel. For instance, for every periodic element in $\Mod(S_7)$, there is an element of $\Mod(S_9)$ of the same order. For Theorems~\ref{thm:finite} and \ref{thm:homeo}, the situation is more extreme, as homeomorphism groups of spheres contain elements of every finite order.

Despite these obstacles, we are still able to apply Theorem~\ref{thm:LM} to these situations by analyzing non-cyclic finite subgroups of $\Mod(S_g)$. Our analysis builds upon prior work on these finite subgroups by a number of authors: May--Zimmerman, M\"uller--Sarkar, Accola, Maclachlan, and Weaver. Our work was facilitated by the cataloging of finite subgroup of $\Mod(S_g)$ up to $g=48$ by Breuer and by Paulhus; see \cite{Breuer, Paulhus, PaulhusWeb}. For Theorems~\ref{thm:finite} and \ref{thm:homeo} we rely on work by a number of authors on classifying finite groups acting on spheres; some is classical and some is relatively recent work by Mecchia, Zimmermann, and Pardon. In working with finite groups we have made use of the computer algebra system GAP, its SmallGroupLibrary, and the web resource GroupNames \cite{GAP, GAPSmallGroups, GroupNames}.

In the remainder of the introduction, we discuss prior results that motivated our work, give some additional context, and point out some questions and conjectures that are related to our results.

\p{Homomorphisms between mapping class groups}

Aramayona--Souto proved a rigidity theorem for homomorphisms between the mapping class groups of surfaces $S_{g,n,b}$, which are surfaces of genus $g$ with $n$ punctures and $b$ boundary components \cite{AS}. The mapping class groups they consider are pure, in that they require mapping classes to fix punctures and boundary components pointwise.

\begin{theorem}[Aramayona--Souto, Theorem~1.1, \cite{AS}]
\label{thm:ASfull}
Let $S~=~S_{g,n,b}$ and $S'~=~S'_{g',n',b'}$, such that $g \geq 6$ and $g' \leq 2g-1$. If $g'=2g-1$, suppose also that $S'$ is not closed. Then every nontrivial homomorphisms $\phi: \Mod(S) \to \Mod(S')$ is induced by an embedding $S \to S'$.
\end{theorem}

Aramayona--Souto also prove that the conclusion of this theorem holds when $g=g' \in \{4,5\}$. They explain the necessity of their upper bound of $2g-1$ by observing that there is a ``double embedding" homomorphism $\Mod(S_{g,0,1}) \to \Mod(S_{2g,0,0})$.

When restricted to the case of closed surfaces, Theorem~\ref{thm:ASfull} says that for $g \geq 6$ and $h < 2g-1$, every homomorphism $\phi: \Mod(S_g) \to \Mod(S_h)$ is trivial, except for the possibility of an isomorphism in the case $g=h$. Some cases of this result for closed surfaces were previously known. The case $h=g$ was previously treated by Ivanov--McCarthy under the further hypothesis that the homomorphism is injective \cite{IvMc}. The range $h<g$ was previously handled by Harvey--Korkmaz \cite{HK}. Our proof strategy for Theorem~\ref{thm:newAS} is similar to their approach: they observe that some non-identity power of elements of order $4g+2$ must lie in the kernel of any map $\phi: \Mod(S_g) \to \Mod(S_h)$, $h<g$; that these powers are all normal generators of $\Mod(S_g)$ except when the power is a hyperelliptic involution; and that an additional argument handles this last case where the map factors through the symplectic representation.

Theorem~\ref{thm:newAS} gives a new proof of Theorem~\ref{thm:ASfull} when restricted to the case where $S$ and $S'$ are closed surfaces and $g \neq g'$. Further, it is an easy consequence of the proof of Theorem~\ref{thm:newAS} that we may replace $S_h$ with a surface of genus $h$ and arbitrarily many punctures and boundary components, since $\Mod(S_{h,n,b})$ has no additional finite subgroups compared to $\Mod(S_h)$. For this replacement, note that it is not necessary for us to restrict to $\Mod(S_{h,n,b})$ being pure. Our Theorem~\ref{thm:newAS} also covers the additional small values of $g$ of 3, 4, and 5, confirming an expectation that Aramayona--Souto state in their paper. We also extend their upper bounds for $g'$ slightly when $g$ is odd. This is possible because, for closed surfaces, their upper bound of $2g-1$ does not have the same natural justification that appears for surfaces with boundary.

It is important to note that, even for closed surfaces, some upper bound is necessary on $h$ to ensure that the homomorphism is trivial. First, since $\Mod(S_g)$ is residually finite, it has a rich supply of finite quotients \cite{Grossman}. As every finite group is a subgroup for some $\Mod(S_h)$, we obtain for all $g>0$ non-trivial homomorphisms $\Mod(S_g) \to \Mod(S_h)$ that factor through finite quotients. An even more striking reason why some upper bound on $h$ is necessary is a result of Aramayona--Leininger--Souto: that for all $g \geq 2$, there exists a nontrivial connected cover $S_h$ of the surface $S_g$ such that $\Mod(S_g)$ injects into $\Mod(S_h)$ \cite{ALS}. All of these sources of nontrivial homomorphisms between mapping class groups are in accord with the conjectural picture proposed by Mirzakhani and recorded in \cite{AS} that every homomorphism between mapping class groups of sufficiently high genus has either finite image or is induced by some manipulation of surfaces. Regarding finite quotients, see also Section 3 of Birman's problem paper \cite{BirmanProblems}.

The preceding results raise the following natural question:

\begin{question}
For each $g \geq 3$, what is the smallest $h>g$ such that there exists a nontrivial homomorphism $\phi:\Mod(S_g) \to \Mod(S_h)$?
\end{question}

We point out that while our Theorem~\ref{thm:newAS} gives approximately the same lower bound on $h$ as that of Aramayona--Souto of about $2g$, our approach suggests that the true bound ought to be higher. Using the work of Breuer and Paulhus that catalogues the finite subgroups of $\Mod(S_g)$ for $g \leq 48$, the following chart indicates for small values of $g$ the smallest $h$ for which $\Mod(S_h)$ contains all of the finite subgroups that are contained in $\Mod(S_g)$. By Theorem~\ref{thm:LM} and our Corollary~\ref{cor:FM}, for all smaller $h$ we have that $\phi:\Mod(S_g) \to \Mod(S_h)$ is trivial as long as $h \leq 3^{g-1}$. This last condition holds for all $g \geq 4$ in the table, while for $\Mod(S_3)$ the table implies that the first candidate target for a non-trivial homomorphism is $\Mod(S_{10})$. We observe that for these limited data points, the first values for $h$ that are candidates for a nontrivial homomorphism $\phi:\Mod(S_g) \to \Mod(S_h)$ are notably larger than $2g$.

\begin{center}
\begin{tabular}{||c | c | c | c | c | c | c||} 
 \hline
 $g$  & 3 & 4& 5& 6& 7& 8 \\ [0.5ex] \hline
 $h$ & 15& 16&21 &$>48$ &$>48$ &40 \\ [0.5ex] 
 \hline
\end{tabular}
\end{center}

\p{Homomorphisms to homeomorphism groups of spheres} 

Franks--Handel prove a number of theorems showing that homomorphisms from mapping class groups to many groups are trivial \cite{FH}. Their main result is that homomorphisms $\Mod(S) \to \GL(n,\mathbb{C})$ are trivial for $S$ a finite-type surface whenever $g \geq 3$ and $n<2g$. They also apply the theorem of Aramayona--Souto to show that for $g\ge 6, 1<g'<2g$, every homomorphism $\Mod(S_g) \to \Homeo(S_{g'})$ is trivial. As these bounds exclude $g'=0,1$, they go on to show that for $g \geq 3$, all homomorphisms from $\Mod(S_g)$ to each of $\Homeo(\Sph^1$) and $\Homeo(S_1)$ are trivial, and they also show the following theorem.

\begin{theorem}[Franks--Handel, Theorem~1.4, \cite{FH}]
\label{thm:FH}
For $g \geq 7$, every homomorphism $\phi: \Mod(S_g) \to \Diff^+(\Sph^2)$ is trivial.
\end{theorem}
Note that Franks–Handel considered the ordinary mapping class groups (rather than the extended mapping class groups) throughout their paper.
Our Theorem~\ref{thm:homeo} includes an extension of this theorem of Franks--Handel to the target $\Homeo^+(\Sph^2)$ and also covers several additional small genus cases. Note that the lower bound of $g \geq 3$ in the statement of Theorem \ref{thm:homeo} is necessary, since in the cases $g=1,2$ there are nontrivial homomorphisms that factor through the abelianization of $\Mod(S_g)$.

Zimmermann showed that every homomorphism $\Mod(S_3) \rightarrow \Diff^+(\Sph^n)$ is trivial when $n \leq 4$ \cite{ZimmermannSp}. His proof approach also involves an analysis of torsion.

Let $\Aut(F_r)$ be the automorphism group of a free group of rank $r$, and let $\SAut(F_r)$ be its unique index 2 subgroup. A corresponding result about homomorphisms $\SAut(F_r) \rightarrow \Homeo^+(\Sph^n)$ was proven by Bridson--Vogtmann \cite{BridsonVogtmann}. They showed that for $n<r+1$, every such homomorphism is trivial. Their result is sharp, as there is a standard linear action of $\SAut(F_r)$ on $\Sph^{r-1}$ that factors through $\SL(r,\R)$. Bridson--Vogtmann apply arguments about torsion to prove their results. One difference between $\Aut(F_r)$ and $\Mod(S_g)$ is that there are natural embeddings $\Aut(F_r)\to \Aut(F_{r+1})$ for all $r\ge 1$, which do not exist for the groups $\Mod(S_g)$. The embedding implies that torsion persists in the former but not the latter, which makes finding obstructions in finite subgroups more difficult for the groups $\Mod(S_g)$. 

In Theorem~\ref{thm:homeo}, one of our target groups is a diffeomorphism group rather than a homeomorphism group. Unlike in the case of $\Homeo^+(\Sph^2)$, for higher-dimensional spheres there exist actions by finite groups that have wildly embedded fixed point sets, and such actions cannot be smooth; see the survey article by Zimmermann for a discussion \cite{ZimSurvey}. Recent work of Pardon has as a consequence that there do not exist any isomorphism types of finite subgroups of $\Homeo^+(\Sph^3)$ that do not also occur in $\Diff^+(\Sph^3)$ \cite{Pardon}; this result allows us to drop the assumption of smoothness in the case of $\Sph^3$. On the other hand, the finite subgroups of $\Homeo^+(\Sph^4)$ are not yet classified. Under the further hypothesis of smoothness, however, there are results classifying finite group actions on $\Sph^4$; we use these in our proof of Theorem~\ref{thm:homeo}. Despite there not being a full classification of the finite subgroups of $\Homeo^+(\Sph^n)$, for some classes of groups there are lower bounds on the $n$ for which these groups act faithfully by homeomorphisms on $\Sph^n$. Results in this direction by Zimmermann are what allow us to prove Theorem~\ref{thm:finite}.

There do exist examples of nontrivial homomorphisms from mapping class groups to homeomorphisms groups of spheres. Again, there exist many that factor through finite quotients. Another example is to take the symplectic representation $\Mod(S_g)\to \Sp(2g, \R)$ given by the action of $\Mod(S_g)$ on the homology of $S_g$. This gives a homomorphism $\Mod(S_g)\to \Homeo^+(S^{2g-1})$, by acting on the oriented projective space of $\mathbb{R}^{2g}$, similar to the standard homomorphism $\SAut(F_n) \rightarrow \Homeo^+(\Sph^{n-1})$ mentioned above. Zimmermann \cite{ZimmermannSp} proved that the minimal $n$ for a nontrivial homomorphism 
$ \Sp(2g, \R)\to  \Homeo^+(\mathbb{S}^{2g-1})$ to exist is when $n=2g-1$.  Another interesting example is that $\Mod(S_g)$ acts by homeomorphisms on $\mathcal{PMF}(S_g) \cong \Sph^{6g-7}$, the space of projective measured foliations on $S_g$. It would be remarkable if these geometric examples are minimal or special in some sense. 

\begin{conjecture}
For $g \geq 3$:\\

\noindent (1) The minimal $n$ such that there exists a nontrivial homomorphism 
\[
 \Mod(S_g)\to  \Homeo^+(\mathbb{S}^n)\] is $n=2g-1$.\\ 
 
\noindent (2) The minimal $n$ such that there exists an injective homomorphism 
\[ \Mod(S_g)\to  \Homeo^+(\mathbb{S}^n)\] is $n=6g-7$. 

\end{conjecture}

Applying some results from number theory, our proof of Theorem~\ref{thm:finite} shows that there exists $c>0$ such that if there is a nontrivial homomorphism $\Mod(S_g)\to  \Homeo^+(\mathbb{S}^n)$ then $n \geq c\sqrt[11]{g}$. See the paragraph after the proof of Theorem~\ref{thm:finite} for some additional details.

\subsection*{Outline} After proving some preliminary lemmas, we prove Theorem~\ref{thm:newAS} in Section~\ref{sec:newAS}. We then prove Theorems \ref{thm:finite} and \ref{thm:homeo} in Section \ref{sec:homeo}. 
 
\subsection*{Acknowledgments} We thank Dan Margalit and Benson Farb for a number of helpful conversations and for comments on a draft of this article. We thank Samuel Taylor for suggesting that we consider maps to homeomorphism groups of higher-dimensional spheres. We thank Mattia Mecchia for a helpful correspondence about these groups and for suggesting the statement and approach to Theorem~\ref{thm:finite}. We also thank J\"{u}rgen M\"{u}ller for suggesting ways to treat some small values of $g$ in the proof of Theorem~\ref{thm:homeo}. Lei Chen acknowledges support from NSF Grant DMS-2005409. Justin Lanier acknowledges support from NSF Grants DGE-1650044 and DMS-2002187.

\section{Homomorphisms between mapping class groups}
\label{sec:newAS}

The strategy for proving Theorem~\ref{thm:newAS} is straightforward. We first show that some periodic element is in the kernel of the given homomorphism $\Mod(S_g) \to \Mod(S_h)$. We then show that this implies that the homomorphism is trivial. We lay the groundwork for these two steps in two lemmas, Lemmas~\ref{lem:MZ}~and~\ref{lem:trivial}. We then prove Theorem~\ref{thm:newAS}.

We first show that a nontrivial periodic element must lie in the kernel of the homomorphism, for the simple reason that there exists a finite subgroup of $\Mod(S_g)$ that does not exist in any $\Mod(S_h)$ for $h$ in the specified range. The finite subgroups that we use lie in two infinite families, one for when $g$ is even, the other for when $g$ is odd. These families of subgroups were studied by May--Zimmerman; we will draw out the salient features of their work and will repeat some of their arguments for the sake of clarity, but we refer the reader to their papers for full details.

For $n \geq 2$, let $DC_n$ be the dicyclic group of order $4n$, given by the presentation
\[
\langle \ x, y \ | \ x^{2n}=1, x^n=y^2, y^{-1}xy=x^{-1} \ \rangle.
\]
May--Zimmerman showed that when $n$ is even, $DC_n$ has \emph{strong symmetric genus} $n$; that is, the first genus $g$ for which $DC_n$ appears as a subgroup of $\Mod(S_g)$ is when $g=n$ \cite[Theorem~1]{MZDC}. Similarly, they show that when $n$ is odd, $C_4 \times D_n$ has strong symmetric genus $n$, where $C_n$ is the cyclic group of order $n$ and $D_n$ is the dihedral group of order $2n$ \cite[Theorem~3]{MZZD}. Let $G$ stand for any one of these finite groups. In proving their results, May--Zimmerman first show that $G$ does in fact appear as a subgroup of the specified mapping class group $\Mod(S_g)$. To guarantee that this is the first appearance, they then give lower bounds on $h$ for any other $\Mod(S_h)$ containing $G$ as a subgroup, showing that $h>g$. In the proof of Lemma~\ref{lem:MZ}, we follow their method and keep track of lower bounds on $h$, showing that there exists a gap between the first appearance of $G$ as a subgroup and its next appearance within the family of mapping class groups.

\begin{lemma}
\label{lem:MZ}
When $g \geq 2$ is even, $DC_g$ appears as a subgroup of $\Mod(S_g)$ and does not appear in any other $\Mod(S_h)$ with $h<2g-1$. When $g \geq 3$ is odd, $C_4 \times D_g$ appears as a subgroup of $\Mod(S_g)$ and does not appear in any other $\Mod(S_h)$ with $h<2g+1$.
\end{lemma}

Before proving Lemma~\ref{lem:MZ}, we require some preliminaries; see the article of Broughton as a reference \cite{Broughton}. Recall that for any faithful orientation-preserving action of a finite group $G$ on a hyperbolic surface $S_g$ by isometries, we have that
\begin{equation}
\label{eq:A}
A=2g_0-2+\sum_{i=1}^r \Big(1-\frac{1}{\lambda_i}\Big)
\end{equation}
where the normalized area $A$ is the hyperbolic area of the quotient orbifold scaled by $\frac{1}{2\pi}$, the $g_0$ is the genus of the quotient orbifold, the $\lambda_i$ are the orders of the cone points, and the $r$ is the number of cone points. The data of this action is called its signature and it is often encoded as $(g_0;\lambda_1,\dots,\lambda_r)$. The normalized area, group order, and genus $g$ of the original surfaces are related by the equation
\[|G| \cdot A =2g-2.\] 
Finally, for a signature of $(g_0;\lambda_1,\dots,\lambda_r)$ to arise from an action of $G$, these values must satisfy the Riemann--Hurwitz equation and there must be elements $a_1,\dots,a_g$, $b_1,\dots,b_g$, $c_1,\dots,c_r \in G$ that together generate $G$, where $|c_i|=\lambda_i$, and that satisfy 

\begin{equation}
\label{eq:product}
\prod_{i=1}^{g} \ [a_i,b_i] \cdot \prod_{j=1}^{r} \ c_j =1.
\end{equation}

\noindent This equation follows from the fact that a finite group acting on a hyperbolic surface must arise as a quotient of a Fuchsian group.

The methods used here to constrain finite group actions on surfaces are similar to those used in the proofs of the classical $84(g-1)$ and $4g+2$ theorems; see, for instance, \cite[Theorems 7.4, 7.5]{FarbMargalit}.



\begin{proof}[Proof of Lemma~\ref{lem:MZ}]
May--Zimmerman argue that when $g \geq 2$ is even and $DC_g$ is a subgroup of $\Mod(S_h)$, either $h=g$ or $h > g$. We will sharpen their second bound to show that in fact either $h=g$ or $h \geq 2g-1$. Their arguments run by showing a dichotomy that either the normalized area $A$ of a fundamental domain of the action of $DC_g$ on $S_h$ satisfies either $A=\frac{1}{2}-\frac{1}{2g}$ or else $A \geq \frac{1}{2}$. The value of $h$ can then be computed from  the Riemann--Hurwitz equation: $h=1+2gA$. We will improve the latter bound to $A \geq \frac{g-1}{g}$, which implies that $h \geq 2g-1$.

We consider cases based on the values of $g_0$ and $r$. By Equation~(\ref{eq:A}), if $g_0 \geq 2$, then $A \geq 2$. Next, if $g_0=1$, then $r \geq 1$ since $g$ is assumed to be at least 2. If $g_0=1$ and $r \geq 2$, then $A \geq 1$. If $g_0=1$ and $r=1$, then $DC_g$ has a generating set $\langle a,b \rangle$ where $[a,b]=c^{-1}$ and $|c|=\lambda_1$.  It is straightforward to check that $|[a,b]|=g$ for any generating pair for $DC_g$. This yields $A \geq \frac{g-1}{g}$.

Finally, consider the case $g_0=0$. Since $A > 0$, we know that $r \geq 3$. If $r \geq 5$, then we have $A \geq \frac{3}{4} \cdot 2 + \frac{1}{2} \cdot 3 - 2 = 1$. Now assume that $r\le 4$. There is a generating set of  $DC_g$ satisfying $\prod_{j=1}^r c_j=1$. Any generating set for $DC_g$ contains at least one generator outside of $\langle x \rangle$, and by the product restriction there are an even number of these in one of our generating sets. Further, each element outside of $\langle x \rangle$ has order 4. We now treat the cases $r=3$ and $r=4$ in turn.

Let $r=3$. Then exactly two $c_i$ lie outside of $\langle x \rangle$ and these have order 4. These must generate $DC_g$, and a short computation shows that their product has order $2g$. Therefore if $r=3$, the signature is $(0;4,4,2g)$, $A=\frac{1}{2}-\frac{1}{2g}$, and $h=g$, as shown by May--Zimmerman.

Finally, let $r=4$. There are either two or four $\lambda_i$ equal to 4 that correspond to elements outside of $\langle x \rangle$. If there are four, then $A=1$. If there are two and the corresponding $c_i$ together form a generating pair for $DC_g$, then the order of their product is relatively prime to $2g$. Then the product of the remaining $c_i$ must equal the inverse of this product and so has the same order. If one of the remaining $\lambda_i$ is either 2 or 3, then the last $\lambda_i$ must be $2g$ or $2g/3$, respectively. (The latter case is only possible when $g$ is a multiple of 3.) These yield $A=\frac{2g-1}{2g}>\frac{g-1}{g}$ and $A=\frac{7}{6}-\frac{3}{2g}$. The latter is at least 1 when $g \geq 12$ and also satisfies the required bound when $g=6$.

The final possibility is that there are exactly two $\lambda_i$ equal to 4 that correspond to elements outside of $\langle x \rangle$ and so that the corresponding $c_i$ together do not form a generating pair for $DC_g$. Again, the product of the remaining $c_i$ must equal the inverse of this product and so have the same order. Since a product in these remaining two $c_i$ must generate $\langle x \rangle$, if one of the remaining $\lambda_i$ is either 2 or 3, then the last $\lambda_i$ must be again be $2g$ or $2g/3$, respectively. The result follows.

\medskip

Similarly, May--Zimmerman argue that when $g \geq 3$ is odd and $C_4 \times D_g$ is a subgroup of $\Mod(S_h)$, either $h=g$ or $h > g$. We will show how their arguments imply that $h=g$ or $h \geq 2g+1$. Their proof proceeds by showing that the normalized area $A$ of a fundamental domain of the action on $S_h$ by $C_4 \times D_g$ is either exactly $\frac{1}{4}-\frac{1}{4g}$ or else $A \geq \frac{1}{2}$, leaving a few of the final cases as an exercise. These bounds imply that either $h=1+4g(\frac{1}{4}-\frac{1}{4g})=g$ or $h \geq 1+4g \cdot \frac{1}{2}=2g+1$ whenever $C_4 \times D_g$ is a subgroup of $\Mod(S_h)$, as desired.

In what follows, we recap the arguments of May--Zimmerman and fill in the cases they leave as an exercise. If $g_0 \geq 1$, then $A \geq \frac{1}{2}$. Further, if $g_0=0$, then $r \geq 3$, and also the $\lambda_i$ are all even. If $r \geq 4$, then $A \geq \frac{1}{2}$. If $r=3$, they argue that if $\lambda_1=2$ then $h=g$, and if $\lambda_1 \geq 6$ then $A \geq \frac{1}{2}$. It remains to treat $\lambda_1 = 4$. For $\lambda_2 \geq 8$, they show that $A \geq \frac{1}{2}$; they leave $\lambda_2=4,6$ as an exercise.

If $\lambda_1=4$ and $\lambda_2=4$, then each generator $c_1,c_2$ must be a generator of $C_4$ times a reflection in $D_n$, and necessarily the second factor of $c_1 \cdot c_2$ must be a rotation of order $n$. But a pair of such elements does not generate $C_4 \times D_g$; instead they generate an index 2 normal subgroup corresponding to a quotient to $C_2$ given by the parity of word length.

Similarly, if $\lambda_1=4$ and $\lambda_2=6$, then $c_1$ must be a generator of $C_4$ times a reflection in $D_n$, and $c_2$ must be either the identity or the square of a generator in the first factor and a rotation of order 3 in the second factor. This last condition is impossible unless $n$ is a multiple of 3. But a rotation that does not generate the rotation subgroup in $D_n$ cannot be a member of a generating pair for $D_n$. So this case is impossible as well.
\end{proof}

We now prepare the second step, showing that a periodic element in the kernel of $\Mod(S_g) \to \Mod(S_h)$ implies that the homomorphism is trivial whenever $h$ is in the specified range, even in the case when it contains a hyperelliptic involution. The following lemma is similar to a result proved and applied by Harvey--Korkmaz to show their result on homomorphisms $\Mod(S_g) \to \Mod(S_h)$ where $g>h$ \cite[Theorem~7]{HK}. Lemma~\ref{lem:trivial} has the advantage of giving a uniform treatment for all $g\geq 3$. 

\begin{lemma}
\label{lem:trivial}
Let $g \geq 3$ and let $\psi: \Sp(S_g) \to G$ be a homomorphism. If $G$ does not contain $(\mathbb{Z}/3\mathbb{Z})^g$ as a subgroup, then $\psi$ is trivial.
\end{lemma}

\begin{figure}[h!]
\centering
\includegraphics[scale=.6]{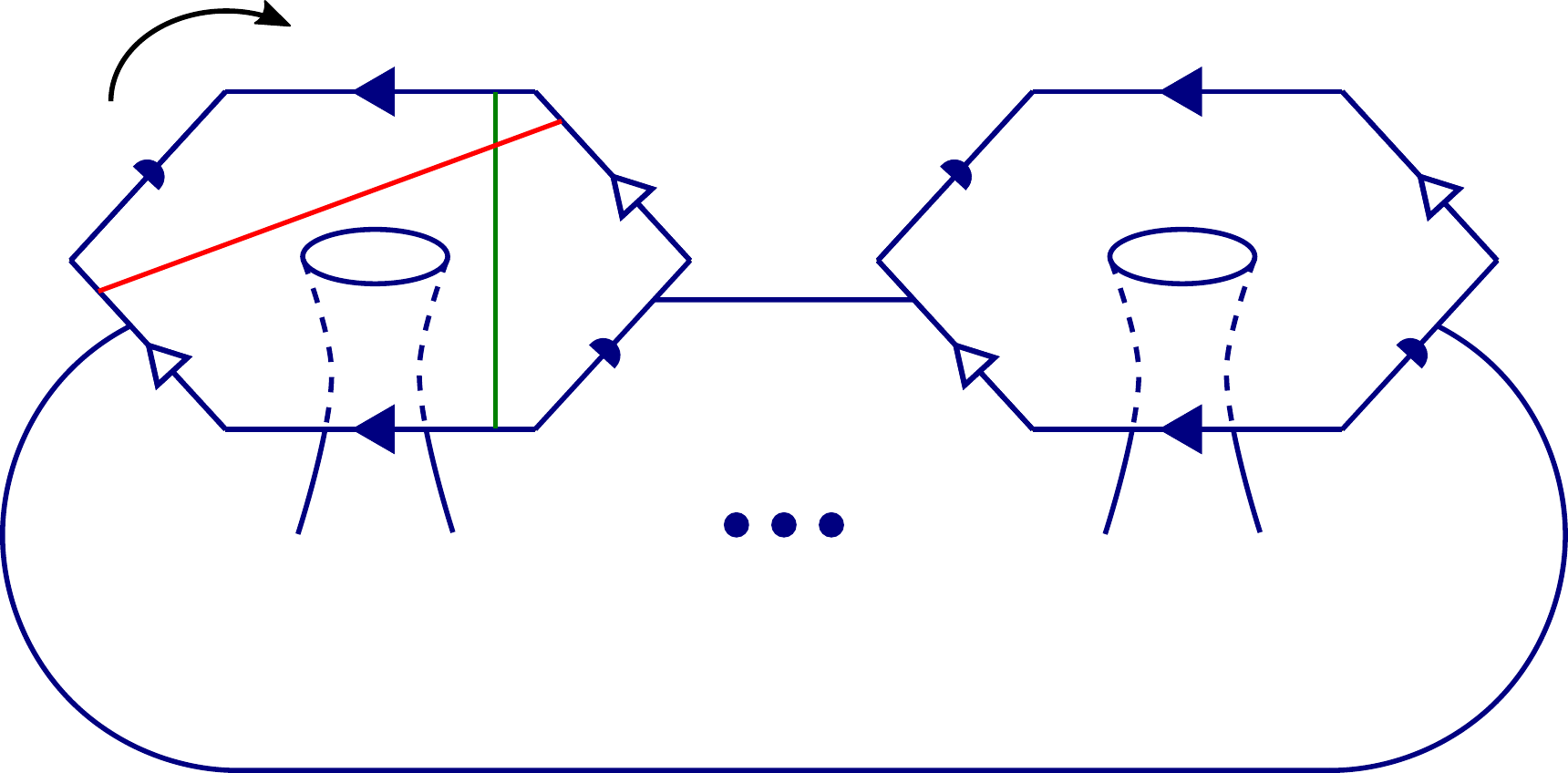}
\caption{The surface $S_g$, where each of $g$ hexagons with sides identified as indicated yields a handle. Each nontrivial element $\tilde{m} \in \tilde{M}$ consists of $1/3$ or $2/3$ rotations about some of the hexagons; we have $i(c,\tilde{m}(c))=1$ for some nonseparating simple closed curve $c$.}
\label{fig:z3}
\end{figure}

\begin{proof}
We proceed by constructing normal generators of $\Sp(2g,\Z)$ and showing that one must be in the kernel of $\psi$.

The group $\Sp(2g,\Z)$ contains $M=(\mathbb{Z}/3\mathbb{Z})^g$ as a subgroup, with generators $m_1$ to $m_g$ of the form 
\[
m_i = \begin{bmatrix} I_{2i-2} & 0 & 0 \\0 & A & 0 \\ 0 & 0 & I_{2g-2i} \end{bmatrix}
\]
Here $A$ is the $2 \times 2$ matrix 
\[
\begin{bmatrix}  -1 & -1 \\ 1 & 0 \end{bmatrix}
\]
of order 3 and $I_n$ denotes the $n \times n$ identity matrix. Every nontrivial element in $M$ is a normal generator of $\Sp(2g,\Z)$, since it is the image of a normal generator of $\Mod(S_g)$, as we now show. The desired elements of $\Mod(S_g)$ are products of roots of Dehn twists about disjoint separating curves, as illustrated in Figure \ref{fig:z3}.
For any nontrivial element $m \in M$, there is a corresponding mapping class $\tilde{m}$ and a nonseparating curve $c$ so that $c$ and $\tilde{m}(c)$ intersect exactly once. By the well-suited curve criterion, it follows that the mapping class $\tilde{m}(c)$  is a normal generator of $\Mod(S_g)$ \cite[Lemma~2.2]{LanierMargalit}; so too is its image $m$ in $\Sp(2g,\Z)$, since normal generators descend to quotients.

As $G$ contains no subgroup isomorphic to $M$, some normal generator of $\Sp(2g,\Z)$ lies in the kernel of $\psi$; therefore $\psi$ is trivial.
\end{proof}

The following corollary can be thought of as a special case of a theorem of Farb--Masur in the case of $\Sp(2g,\Z)$; they show that for any irreducible lattice $\Gamma$ in a semisimple lie group $G$ of $\R$-rank at least two, the image of any homomorphism $\phi: \Gamma \to \Mod(S)$ is finite \cite[Theorem~1.1]{FarbMasur}.

\begin{corollary}
\label{cor:FM}
Let $g \geq 3$ and let $h \leq 3^{g-1}$. Then every homomorphism $\psi: \Sp(2g, \Z) \to \Mod(S_h)$ is trivial.
\end{corollary}

\begin{proof}
This follows immediately from Lemma~\ref{lem:trivial} by a result of M\"uller--Sarkar, who showed that the strong symmetric genus of $(\mathbb{Z}/3\mathbb{Z})^g$ is $1+3^{g-1} \cdot \mu_0(g)$, where $\mu_0(g) \geq 1$ when $g \geq 3$ \cite[Section~9.1]{MullerSarkar}.
\end{proof}

With these preliminaries established, the proof of Theorem~\ref{thm:newAS} is straightforward.

\begin{proof}[Proof of Theorem~\ref{thm:newAS}]
Let $g \geq 3$ be even and $g<h<2g-1$. By Lemma~\ref{lem:MZ}, $\Mod(S_g)$ contains $DC_g$ as a subgroup and $\Mod(S_h)$ does not. Then for any homomorphism $\phi: \Mod(S_g) \to \Mod(S_h)$, a nontrivial periodic element $f$ lies in the kernel. If $f$ is not a hyperelliptic involution, we conclude that $\phi$ is trivial by Theorem~\ref{thm:LM}. If $f$ is a hyperelliptic involution, then $\phi$ factors through the symplectic representation. Since $h<2g-1 \leq 3^{g-1}$, we conclude that $\phi$ is trivial by Corollary~\ref{cor:FM}.

The proof for odd $g \geq 3$ proceeds in the same way using the finite group $C_4 \times D_g$.
\end{proof}

\section{Homomorphisms to homeomorphism groups of spheres}
\label{sec:homeo}
In this section we prove Theorems~\ref{thm:finite} and \ref{thm:homeo}, which constrain homomorphisms from mapping class groups to homeomorphism groups of spheres.

Let $\mathbb{Z}_k$ be the cyclic group of order $k$. Zimmermann showed that there are strong restrictions on which spheres certain metacyclic groups $\mathbb{Z}_p\rtimes\mathbb{Z}_k$ have a faithful action, where $p$ is an odd prime \cite{Allsphere}. While Zimmermann throughout his paper assumes smoothness of actions, the results that are relevant for us hold true for general actions by homeomorphisms. We summarize these here. Zimmermann shows the following preliminary homological proposition about actions of $\mathbb{Z}_p\rtimes\mathbb{Z}_k$ on homology spheres using a spectral sequence argument.

\begin{proposition}[Zimmermann, Proposition~2, \cite{Allsphere}]
\label{prop:Zprelim}
Let $p$ be an odd prime and $k$ be a positive integer. Let $G=\mathbb{Z}_p\rtimes\mathbb{Z}_k$ be a semidirect product with an effective action of $\mathbb{Z}_k$ on the normal subgroup $\mathbb{Z}_p$. Suppose that $G$ admits a faithful action on a manifold $M$ with the mod $p$ homology of $\Sph^n$ and that the group $\mathbb{Z}_p$ acts freely on $M$. Then $n+1$ is a multiple of $2k$ if all elements of $G$ act as the identity on $H^n(M;\mathbb{Z}_p)\cong \mathbb{Z}_p$ (the orientation-preserving case), or an odd multiple of $k$ if some element of $G$ acts as the minus identity (the orientation-reversing case).
\end{proposition}

We now give a proof of one of Zimmermann's main results \cite[Theorem 1(i)]{Allsphere}, adapted to apply to actions by general homeomorphisms.

\begin{theorem}
\label{thm:pqbound}
For $p$ an odd prime and $k$ an integer, let $G=\mathbb{Z}_p\rtimes\mathbb{Z}_k$ be a semidirect product with an effective action of $\mathbb{Z}_k$ on the normal subgroup $\mathbb{Z}_p$. 
\begin{enumerate}
\item If $k$ is odd, the minimal $n$ such that there is a faithful topological action of the group $G$ on $\mathbb{S}^n$ is $n= 2k-1$.
\item If $k$ is even, the minimal $n$ such that there is a faithful orientation-preserving topological action of the group $G$ on $\mathbb{S}^n$ is $n= k$. 
\item
If $k$ is even, the minimal $n$ such that there is a faithful topological action of the group $G$ on $\mathbb{S}^n$ is $n= k-1$. 
\end{enumerate}
\end{theorem}

\begin{proof}
The group $G$ contains a normal subgroup $Z\cong \mathbb{Z}_p$. Let $F$ be the fixed point set of the action of $Z$ on $\Sph^n$. If $F$ is empty, then Proposition~\ref{prop:Zprelim} implies that $n+1$ is a multiple of $2k$, and we have $n+1 \ge 2k$. If $F$ is not empty, then $F$ is a $\mathbb{Z}_p$ homology sphere by Smith theory. By the Alexander Duality theorem, its complement in $\Sph^n$ has the same mod $p$ homology as a sphere of dimension less than $n$. Since $Z$ acts on $\mathbb{S}^n-F$ freely and $\mathbb{S}^n-F$ is $G$-invariant, again Proposition~\ref{prop:Zprelim} implies that $n+1 \ge 2k$ or $n+1\ge k$ depending on whether the action on the top homology is trivial or not; only the former may occur when $k$ is odd. Finally, $G$ admits a faithful linear action on $\Sph^{2k-1}$ if $k$ odd and on $\Sph^{k-1}$ and $\Sph^k$ if $k$ is even, the latter if the action is required to preserve orientation; see for instance \cite[Example 9.2.3, p.155]{Brown}.\end{proof}

We are now prepared to prove Theorem \ref{thm:finite}.

\begin{proof}[Proof of Theorem~\ref{thm:finite}]
For any $n$ and a prime number $q>\frac{n+1}{2}$, there exists a prime number $p=kq+1$ for some natural number $k$ by Dirichlet's theorem on arithmetic progressions (see, e.g., \cite{Dirichlet} for an analytic proof). Therefore, $\mathbb{Z}_q$ admits an effective action on $\mathbb{Z}_p$ since $\mathbb{Z}_q <\text{Aut}(\mathbb{Z}_p)=\mathbb{Z}_{p-1}$. Let $G=\mathbb{Z}_p\rtimes\mathbb{Z}_q$ be a semidirect product with an effective action of $\mathbb{Z}_q$ on the normal subgroup $\mathbb{Z}_p$.  Then by Theorem~\ref{thm:pqbound}, we know that $G$ is not a subgroup of $\Homeo^+(\mathbb{S}^n)$. As observed by Weaver \cite{Weaver}, a result of Kulkarni \cite{KulkarniArith} implies $G$ is a subgroup of $\Mod(S_g)$ for all sufficiently large $g$. Then for all sufficiently large $g$ and for any $\phi: \Mod(S_g) \rightarrow \Homeo^+(\Sph^n)$ there is a nontrivial periodic element in its kernel. Every nontrivial element of $G$ has odd order, so Theorem~\ref{thm:LM} implies that $\phi$ is trivial.
\end{proof}

Following the steps of the proof of Theorem~\ref{thm:finite} yields the following rough bound on the biggest $g$ such that $\Mod(S_g)$ can act nontrivially on $\Sph^n$ by homeomorphisms. Given $n$, Bertrand's postulate yields a prime number $\frac{n+1}{2}<q<n+1$. Linnik's theorem and its generalizations give bounds on the smallest prime $p\equiv 1 \mod q$ namely there exist $c$ and $L$ such that there exists such a prime number $p\equiv 1 \mod q$ such that $p < cq^L$. The best known bound on $L$ currently known is 5, due to Xylouris \cite{Linnik}. (Better bounds on $L$ are known under the assumption of the Generalized Riemann Hypothesis \cite{GRH}.) Finally, Weaver \cite{Weaver} computed the stable upper genus of $\mathbb{Z}_p\rtimes\mathbb{Z}_q$; the dominating term in the resulting formula is $p^2q$. Therefore the largest $g$ such that $\Mod(S_g)$ has a nontrivial action on $\Sph^n$ by homeomorphisms satisfies \[
g = O(p^2q) = O((q^5)^2q) = O(q^{11}) = O(n^{11}).\]

For some small values of $n$, we are able to give sharp or nearly-sharp bounds on the values of $g$ in Theorem~\ref{thm:finite} that yield trivial homomorphisms. In proving Theorem~\ref{thm:homeo}, we consider each of the four target groups in turn: $\Homeo^+(\Sph^2)$, $\Homeo^+(\R^3)$, $\Homeo^+(\Sph^3)$, and then $\Diff^+(\Sph^4)$ and $\Homeo^+(\Sph^4)$. The constraining finite subgroups are the same for the first two target groups; for the latter target groups, a deeper analysis of finite subgroups is required.

We begin by introducing a class of groups that will play an important role in our proofs. The split metacyclic group $D_{p,q}$ of order $pq$ with $p$ and $q$ two prime numbers is the group with presentation
\begin{align}
D_{p,q} = \langle a,b \ | \ a^p=b^q=1, bab^{-1}=a^r \rangle=\mathbb{Z}_p\rtimes \mathbb{Z}_q
\end{align}
where $r$ is a solution (other than 1) to the congruence $r^q \equiv 1$ (mod $p$). Such a solution exists exactly when $p$ is 1 mod $q$, and different solutions yield isomorphic groups. Note that $D_{p,2}$ is a dihedral group. 

\begin{proof}[Proof of Theorem \ref{thm:homeo}] Let $g \geq 3$ and let $\phi: \Mod(S_g) \rightarrow \Homeo^+(\Sph^2)$ be a homomorphism. It is a classical result of Brouwer, Eilenberg, and de Ker\'ekj\'art\'o \cite{Brouwer,K,Eilenberg} that every finite subgroup of $\Homeo^+(\Sph^2)$ is conjugate to a finite subgroup of $\SO(3)$. These are the cyclic groups $C_n$, the dihedral groups $D_n$, and the tetrahedral, octahedral, and icosahedral groups $A_4$, $\Sigma_4$, and $A_5$. This classification goes back to the work of Klein; see \cite[Chapter 19]{Burnside} or \cite[Sect.I.3.4]{Breuer} for treatments.

When $g \geq 3$, $\Mod(S_g)$ contains a finite subgroup that is not isomorphic to any subgroup of $\SO(3)$. For instance, we have the Accola--Maclachlan subgroup $C_2 \times C_{2g+2}$ of $\Mod(S_g)$, which attains the $4g+4$ bound on its largest abelian subgroup (see, for instance, \cite{WeaverSurvey}). This implies that $\phi$ has a nontrivial periodic element in its kernel. If this element is not a hyperelliptic involution, Theorem~\ref{thm:LM} implies that $\phi$ is trivial. Otherwise, $\phi$ factors through the symplectic representation. Since $\SO(3)$ does not contain $(C_3)^g$ as a subgroup for $g \geq 3$, Lemma~\ref{lem:trivial} implies that $\phi$ is again trivial.
\medskip

We next consider homomorphisms to $\Homeo^+(\R^3)$. Zimmermann showed that every orientation-preserving action by homeomorphisms of a finite group on $\Sph^3$ that has a global fixed point is a finite subgroup of $\SO(3)$; he therefore also concludes that every finite subgroup of $\Homeo^+(\R^3)$ is a finite subgroup of $\SO(3)$ \cite[Corollary~1]{ZimS3}. This was also shown independently by Kwasik--Sun \cite{KwasikSun}. Therefore we may conclude just as we did for $\Homeo^+(\Sph^2)$ above that $\phi$ is trivial. 
\medskip

Next we consider homomorphisms to $\Homeo^+(\Sph^3)$. By recent work of Pardon, every continuous action of a finite group on a smooth 3-manifold is a uniform limit of smooth actions \cite{Pardon}. Pardon's result guarantees that the isomorphism types of the finite subgroups of $\Homeo^+(\Sph^3)$ are identical to those of $\Diff^+(\Sph^3)$. By the Geometrization Theorem, every finite group acting smoothly or locally linearly on $\Sph^3$ is geometric, that is, it is conjugate to a finite subgroup of $\SO(4)$. Every finite subgroup of $\SO(4)$ is a subgroup of some central product $P_1 \times_{C_2} P_2$, where each $P_i$ is one of the binary polyhedral groups: the cyclic groups $C_{2n}$, the binary dihedral groups $D_n^*$, and the binary tetrahedral, binary octahedral, and binary icosahedral groups $A_4^*$, $\Sigma_4^*$, and $A_5^*$. These facts are exposited in a survey article by Zimmermann \cite{ZimSurvey}, and a list of subgroups of the binary groups is in the Appendix of \cite{BinaryPolyBook}. (Note that $D_n^* \cong DC_n.)$

We must therefore produce for each $g \geq 3$ a finite subgroup $G_g$ of $\Mod(S_g)$ that is not a subgroup of any $P_1 \times_{C_2} P_2$. (Note that the Accola--Maclachlan subgroup $C_2 \times C_{2g+2}$ no longer suffices.) Producing such subgroups proves the theorem, for then either a nontrivial nonhyperelliptic periodic element is in the kernel of $\phi$, so that $\phi$ is trivial; otherwise $\phi$ factors through the symplectic representation, and since the group $(C_3)^g$ for $g \geq 3$ is not a subgroup of any $P_1 \times_{C_2} P_2$, we conclude by Lemma~\ref{lem:trivial} that $\phi$ is trivial. Indeed, apart from $(\Z/2\Z)^3$, the finite abelian subgroups of $\SO(4)$ can be written as the direct product of at most two finite cyclic groups. 

For $g \geq 42$, we make take $G_g$ to be the split metacyclic group $D_{7,3}$. Weaver computed the stable upper genus of all split metacyclic groups, and in particular he showed that $\Mod(S_g)$ contains $D_{7,3}$ for all $g \geq 42$ \cite[Corollary~4.8]{Weaver}. Weaver additionally shows how to compute the genus spectra of all split metacyclic groups \cite[Theorem~4.7]{Weaver}. Using his formulas, it is straightforward to compute that $\Mod(S_g)$ contains at least one of $D_{7,3}$, $D_{13,3}$, $D_{19,3}$, $D_{31,3}$, or $D_{37,3}$ as a subgroup for all $g \geq 3$ except for the following ten values:
\begin{align}
\{ {4, 5, 7, 11, 13, 16, 23, 25, 34, 41}\}
\end{align}

\noindent For small values of $g$ this data also appears in the catalog of Breuer and Paulhus \cite{Breuer, Paulhus, PaulhusWeb}. On the other hand, $D_{p,3}$ is a not a subgroup of $\Homeo^+(\Sph^3)$ by Theorem~\ref{thm:pqbound}, since $n=3<5=2 \cdot 3-1$. Alternatively, we may see this because all finite subgroups of $\SO(4)$ of odd order are abelian. (The classification of finite subgroups of $\SO(4)$ goes back to work of Seifert--Threlfall \cite{TS1,TS2}; see the paper of Mecchia--Seppia for a contemporary treatment \cite{MeccSepp}.) We have therefore found the required $G_g$ for all but these ten values of $g$. For $g = 4, 11, 16, 25, 34,$ and $41$, we may take $\ G_g = C_5 \rtimes C_4 = \text{SmallGroup($20$,$3$)}$; for these values of $g$ this finite group is a subgroup of $\Mod(S_g)$, but it is not a subgroup of $\Homeo^+(\Sph^3)$ by Theorem~\ref{thm:pqbound}. The latter description of these groups is from the SmallGroupLibrary in the computer algebra system GAP \cite{GAP, GAPSmallGroups}. 

For the remaining four values of $g$, we may proceed as follows, as suggested to us by J\"{u}rgen M\"{u}ller. We begin by setting $G_g =  C_4 \rtimes C_4 = \text{SmallGroup($16$,$4$)}$ for $g = 5, 7, 13,$ and $23$. We have that $G_g$ is a subgroup of $\Mod(S_g)$ by consulting the catalogue of finite subgroups of $\Mod(S_g)$ up to $g=48$ by Breuer. To see that the group $C_4 \rtimes C_4$ is not a subgroup of $\SO(4)$, it suffices to check that they have no real faithful characters in dimension at most 4, realizable over $\R$ and with trivial determinant character. This is a finite check using the groups' character tables and where realizability can be determined using the Frobenius-Schur indicator. Note that it is only for these four values of $g$ where the argument invokes the result of Pardon.

\medskip

Finally, let $\phi: \Mod(S_g) \rightarrow \Diff^+(\Sph^4)$ be a homomorphism. A theorem of Mecchia--Zimmermann \cite{MeZi} states that every finite group that acts smoothly on $\Sph^4$ and preserves orientation lies on the following list:

\begin{enumerate}
\item orientation-preserving finite subgroups of $\Orth(3) \times \Orth(2)$ and of $\Orth(4) \times \Orth(1)$,
\item orientation-preserving subgroups of the Weyl group $W=(C_2)^5 \rtimes \Sigma_5$,
\item $A_5$, $\Sigma_5$, $A_6$, $\Sigma_6$, and
\item finite subgroups of $\SO(4)$ and 2-fold extensions thereof.
\end{enumerate}

Any group satisfying these conditions either is a subgroup of $\SO(5)$ or has a subgroup of index 2 in $\SO(4)$. Since the group $(C_3)^g$ is not on this list for $g \geq 3$, as with target $\Homeo^+(\Sph^3)$ it again suffices to produce for $g \geq 5$ a $G_g$ that is a subgroup of $\Mod(S_g)$ and that is not a subgroup of $\Diff^+(\Sph^4)$. For all but the ten exceptional values of $g$ we may take $G_g$ to be some $D_{p,3}$, just as in the previous case. We have already established that they are subgroups of $\Mod(S_g)$. To see that these $G_g$ are not subgroups of $\Diff^+(\Sph^4)$, nor in fact of $\Homeo^+(\Sph^4)$, it is again enough to appeal to Theorem~\ref{thm:pqbound}. Here we do not require the claim about $(C_3)^g$, since these finite subgroups have odd order. It follows that for $g \geq 42$, every homomorphism $\Mod(S_g) \rightarrow \Homeo^+(\Sph^4)$ is trivial.

For the remaining cases, we may take the following, again as suggested by J\"{u}rgen M\"{u}ller:
\begin{align*}
\text{for } g = 5, 13, 25, 41:& \ G_g = (C_4 \times C_2) \rtimes C_4 = \text{SmallGroup($32$,$2$)}\\
\text{for } g = 7, 16, 34:& \ G_g = C_9 \rtimes C_3       = \text{SmallGroup($27$,$4$)}\\
\text{for } g = 23:& \ G_g = \SL_2(3) \rtimes C_4   = \text{SmallGroup($96$,$66$)}
\end{align*}

\noindent Again, we can appeal to Breuer's catalog to see that $G_g$ is a subgroup of $\Mod(S_g)$, and the same character theory argument shows that these groups are not subgroups of $\SO(5)$, nor do they have a subgroup of index 2 in $\SO(4)$. By the theorem of Mecchia--Zimmermann, we then have than these $G_g$ do not lie in $\Diff^+(\Sph^4)$.

Finally, for $g=11$, every finite subgroup of $\Mod(S_g)$ lies on the list of Mecchia--Zimmermann. However, using their approach and their preliminary results, we now show that the following group does not lie in $\Homeo^+(\Sph^4)$:
\[G=G_{11}= C_5 \rtimes C_8=\text{SmallGroup($40$,$3$)}\]

\noindent Set $a$ and $b$ as generators of the two factors of $C_5 \rtimes C_8$, respectively. A finite $p$-group, for $p$ an odd prime, acting on a sphere has as its fixed point set a sphere of even codimension, while for $p=2$ this holds for cyclic $2$-groups acting orientation-preservingly and freely outside the fixed point set; for instance, see \cite[Sect.IV.4.4, 4.5, 4.7]{Borel}. If the action of $G$ is faithful, then we have that $\Fix \langle a \rangle$ is either an $\Sph^2$ or an $\Sph^0$.

In the first case,  by \cite[Lemma~3]{MeZi}, we know that $G < \text{S}(\text{O}(3) \times \text{O}(2)) < \SO(5)$, but the same character theory shows that $G$ is not in $\SO(5)$, a contradiction.

In the second case, $\Fix \langle a \rangle$ is an $\Sph^0$. Let $H\le G$ be the pointwise stabilizer of $\Fix \langle a \rangle$, which is at most index $2$ in $G$.  By \cite[Lemma~2]{MeZi} we have that $H < \SO(4)$. However, by the same character theory, we know that $G$  is not a subgroup of $\SO(4)$. Therefore $H$ is the unique index $2$ subgroup of $G$, which is isomorphic to $C_5\rtimes_2 C_4 = \text{SmallGroup($20$,$1$)}$ where the action of $C_4$ on $C_5$ is through $C_2$.

Now we analyze $\Fix\langle b^4 \rangle$, which contains $\Fix\langle a \rangle$ because $b^4\in H$. If  $\Fix \langle b^4 \rangle$ is an $\Sph^2$, then by \cite[Lemma~3]{MeZi}, we know that $G < \SO(5)$, then the same character theory shows that $G$ is not in $\SO(5)$, a contradiction.

Therefore $\Fix\langle b^4 \rangle$ is an $\Sph^0$, which implies that $\Fix\langle b^4 \rangle=\Fix\langle a \rangle$. We know that $\Fix\langle b \rangle$ is subset of $\Fix\langle b^4 \rangle$ by the fact $b$ and $b^4$ commute and $\Fix\langle b \rangle$ is not empty by Lefschetz fixed point Theorem. Therefore \[
\Fix\langle b \rangle=\Fix\langle b^4 \rangle=\Fix\langle a \rangle,\]
which contradicts to the fact that $H$, the stabilizer of $\Fix\langle a \rangle$, is a proper subgroup of $G$. Therefore $G$ is not a subgroup of $\Diff^+(\Sph^4)$, and so every homomorphism $\Mod(S_{11}) \rightarrow \Diff^+(S^4)$ is trivial. 
\qedhere

\end{proof}

\bibliographystyle{abbrv}
\bibliography{constrain}

\end{document}